\documentclass[reqno, 12pt]{amsart}
\usepackage {amsmath}
\usepackage{amsthm, amscd, amsfonts, amssymb, graphicx, color}
\usepackage[bookmarksnumbered, colorlinks, plainpages]{hyperref}

\newtheorem{lem}{\textbf{Lemma}}[section]
\newtheorem{thm}[lem]{\textbf{Theorem}}

\newtheorem{rem}[lem]{\textbf{Remark}}
\newtheorem{com}[lem]{\textbf{Comment}}

\theoremstyle{definition}

\theoremstyle{definition}
\theoremstyle{remark}

\begin{document}
\title[A simple     proof for     Kazmi et al.'s  iterative
scheme]{A simple     proof for     Kazmi et al.'s  iterative
scheme}
\author{E. Soori and  Ravi P. Agarwal}
\address[affil]{Department  of Mathematics, Lorestan University, P.O. Box 465, Khoramabad, Lorestan, Iran(Soori)\\
Department of Mathematics, Texas A and M University-Kingsville, Texas 78363, USA(Agarwal)}
\email{sori.e@lu.ac.ir, Ravi.Agarwal@tamuk.edu}

\date{}

\begin{abstract}
In this paper we provide a simple proof of the existence   iterative
scheme  by using   two Hilbert spaces   due to Kazmi et al.   [K. R. Kazmi, R. Ali, M. Furkan, Hybrid iterative method for split monotone \ldots, Numer Algor, 2017].
\end{abstract}


\maketitle


\section{\textbf{Introduction and methods  }}
In this section, we indicate the wrong part of the proof of theorem 3.1 of \cite{kaz}.
In   following comment, we refer the readers to the line 8 in the page 15 in the below of the equation (3.41) of \cite{kaz}. We claim that this part of the proof is wrong, then we modify the proof.

\begin{thm} \cite[Theorem 3.1]{kaz}
   Let $H_{1}$ and $H_{2}$ be real Hilbert spaces and $C \subseteq H_{1}$ and $Q \subseteq H_{2}$ be
nonempty, closed, and convex sets. Let $A : H_{1}\rightarrow H_{2}$ be a bounded linear operator
with its adjoint operator $A^{*}$;  let $M_{1} : H_{1} \rightarrow 2^{H_{1}}$ , $M_{2} : H_{2} \rightarrow 2^{H_{2}}$ be multi-valued
maximal monotone operators. Let $f : C \rightarrow H_{1}$ and $g : Q \rightarrow H_{2}$ be $\theta_{1}$- and $\theta_{2}$-inverse strongly monotone mappings, respectively, let $S : C \rightarrow H_{1}$ be a nonexpansive
mapping and let $\{T_{i} \}^{\mathbb{N}}_{i=0}: C \rightarrow C$ be a finite family of nonexpansive mappings. For
every $n \in \mathbb{N}\cup \{0\}$, let $W_{n}$ be a$W$-mapping generated by $T_{1}, \cdots, T_{N}$ and $\lambda_{n,1}, \cdots, \lambda_{n,N} $.
Assume that 	$\Gamma= \Omega \bigcap \Phi \neq \emptyset$. Let the iterative sequences $\{u_{n}\}$, $\{y_{n}\}$ and $\{x_{n}\}$ be
generated by hybrid iterative algorithm:
\begin{align}
   x_{0} \in C, \quad C_{0} = C   \nonumber   &\\
  u_{n} = (1 - \alpha_{n})x_{n} + \alpha _{n}P_{C}(\sigma _{n} Sx_{n} + (1 - \sigma_{n})W_{n}x_{n}); \nonumber & \\
   z_{n} = U(u_{n}); \quad w_{n} = V (Az_{n});\quad y_{n} = z_{n} + \gamma A^{*}(w_{n} - Az_{n}); \nonumber & \\
  C_{n} = \{z \in C :  \|y_{n} - z \|^{2} \leq (1 - \alpha_{ n}\sigma _{n})\| x_{n} - z \|^{2} + \alpha_{n}\sigma_{n}\| Sx_{n} - z\| ^{2}\}; \nonumber & \\
  Q_{n} = \{z \in C :  \langle x_{n} - z, x_{0} - x_{n}  \rangle \geq0\}; \nonumber &    \\
  x_{n+1} = P_{C_{n}\bigcap Q_{n}}x_{0}, n \geq 0.
\end{align}

where $U := J^{M_{1}}_{\lambda}  (I -\lambda f )$, $V := J^{M_{2}}_{\lambda} (I -\lambda g)$, $A(Range(U)) \subseteq  Q$, $\gamma \in (0, \frac{1}{\|A\|^{2}}) $.

Let $\{\lambda_{n,i} \}_{i=1}^{N}$
  be a sequence in $[0, 1]$ such that $\lambda_{n,i} \rightarrow \lambda_{i}$, $(i = 1, 2, \cdots,\mathbb{N})$ and $\lambda \in (0, \alpha)$ with $\alpha = 2\min\{\theta_{1}, \theta_{2}\}$. Let $\{\alpha_{n}\}$, $\{\sigma_{n}\}$ be real sequences in $(0, 1)$ satisfying
the conditions:
\begin{enumerate}
  \item [(i)] $\displaystyle\lim _{n\rightarrow \infty} \sigma_{n} = 0$,
  \item [(ii)] $\displaystyle\lim_{n\rightarrow \infty}  \frac{\|x_{n}- u_{n}\|}{\alpha_{n} \sigma _{n}}=0$.
\end{enumerate}

Then $\{x_{n}\}$ converges strongly to $z \in \Gamma $, where $z = P_{\Gamma} x_{0}$.
\end{thm}

   \begin{rem}
    The authors claim the weak convergence \eqref{fgg}. Moreover, the content of the  line 3 in the page 16 in the above  of the equation (3.43) is wrong, because  the authors take limit on $n$ to get the equation (3.43) while  they have  used  from the strong  convergence  $\frac{x^{*} - u_{n}-x_{n}}{\alpha_{n}}-x_{n}   $ instead of the weak  convergence of that. Note that, if   $x^{*} -\frac{ u_{n}-x_{n}}{\alpha_{n}}-x_{n}   $ converges strongly to $0$, since moreover,    $\displaystyle \lim_{n\rightarrow \infty} \frac{\| u_{n}-x_{n}\|}{\alpha_{n}}=0$  (the line 6 in the page 15),  then we have
 \begin{equation*}
   \lim_{n\rightarrow \infty} \|x^{*}-x_{n}\| \leq  \lim_{n\rightarrow \infty} \|x^{*} -\frac{ u_{n}-x_{n}}{\alpha_{n}}-x_{n} \|+  \lim_{n\rightarrow \infty}\frac{ \|u_{n}-x_{n}\|}{\alpha_{n}}=0,
 \end{equation*}
 then $\{x_{n}\}$ converges strongly to $x^{*}$ while the aim of theorem 3.1 is to prove that $\{x_{n}\}$ converges strongly to $x^{*}$ in the step VI. Hence this is a scientific error in the article.
 \end{rem}
 \begin{rem}
 First, note that in  the line 8 in the page 15 of \cite{kaz}, the authors  claim that
  \begin{equation}\label{fgg}
    \frac{\| u_{n}-x_{n}\|}{\alpha_{n}}+x_{n} \rightharpoonup x^{*}
  \end{equation}
  but this is not valid since we can't add a real  number  with a member of a Hilbert space in general
   \end{rem}
\section{\textbf{Results} }
 Now in this section, we modify the proof  of theorem 3.1 of \cite{kaz}, in the following comment.
 \begin{com}
 First, note that in line 8 of page 15,  instead of the conclusion     $\frac{\| u_{n}-x_{n}\|}{\alpha_{n}}+x_{n} \rightharpoonup x^{*}$, we may  put
 $\frac{  u_{n}-x_{n}}{\alpha_{n}}+x_{n} \rightharpoonup x^{*}$. Indeed, from the fact that  $x_{n} \rightharpoonup x^{*}$ (line 10 page 14) and  $\displaystyle \lim_{n \rightarrow \infty}\frac{\| u_{n}-x_{n}\|}{\alpha_{n}}=0$ (line 6 page 15), we have
 \begin{align*}
 \lim_{n \rightarrow \infty} \langle \frac{  u_{n}-x_{n}}{\alpha_{n}}+x_{n} -x^{*}, y \rangle &= \lim_{n \rightarrow \infty} \langle \frac{  u_{n}-x_{n}}{\alpha_{n}} , y \rangle + \lim_{n \rightarrow \infty} \langle x_{n} -x^{*}, y \rangle \\&\leq \lim_{n \rightarrow \infty} \frac{ \|u_{n}-x_{n}\|}{\alpha_{n}}\|y\|  + \lim_{n \rightarrow \infty} \langle x_{n} -x^{*}, y \rangle=0
 \end{align*}
for each $y \in H$ ($H$ is a real Hilbert space). Now, we should change the line 3 in page 16 in the above of the equation (3.43) by the following:
\begin{equation*}
\leq \|Sx_{n} - Sx^{*}\|\|x- \frac{  u_{n}-x_{n}}{\alpha_{n}}-x_{n}\|+  |\langle Sx^{*} , x^{*}-\frac{  u_{n}-x_{n}}{\alpha_{n}}-x_{n}  \rangle|.
\end{equation*}
 \end{com}
 There is another error  that we mention in the following comment:
\begin{com}
  Another problem, is that the authors have used from the continuity of    $S$   from the weak topology to the norm topology in line 4 of page 16, while this can not be valid in general, then we should to consider  $S$   as a continuous mapping from the weak topology to the norm topology in theorem 3.1.
\end{com}
 \section{\textbf{Discussion}}
In this paper, we find  some  scientific and substantial errors in  \cite{kaz} and then we modify them.   Unfortunately,  the  main result of the  published version of the article is not valid. Then in this short note we help to the paper to receive to    the virtual version.

\section{\textbf{Conclusion}}
In this paper, we find  some  scientific   errors in  \cite{kaz} and then we modify them to rise the virtual   quality of the paper.

\section{ \textbf{Abbreviations}}
Not applicable
\section{\textbf{Declarations} }
\subsection{Availability of data and material}
Please contact author for data requests.
\subsection{Funding}
Not applicable
\section{  \textbf{Acknowledgements}}
 The  author is  grateful to  the University of Lorestan for their support.
 \section{\textbf{Competing interests}}
 The author declares that he has no competing interests.
 \section{\textbf{Authors' contributions}}
 All authors contributed equally to the manuscript, read and approved the final manuscript.

\end{document}